\documentclass[12pt]{amsart}
\usepackage{amsmath} \usepackage{amssymb} \usepackage{amsthm}\usepackage{amscd}
\newtheorem{thm}{Theorem}[section]
\newtheorem{lemma}[thm]{Lemma}
\newtheorem{prop}[thm]{Proposition}
\newtheorem{cor}[thm]{Corollary}

\newtheorem{introthm}{Theorem}

\theoremstyle{definition}
\newtheorem{defn}[thm]{Definition}
\newtheorem{exam}[thm]{Example}

\newtheorem{quest}[thm]{Question}
\newtheorem{introquest}[introthm]{Question}
\theoremstyle{remark}
\newtheorem{remark}[thm]{Remark}
\numberwithin{equation}{section}

\DeclareMathOperator{\Assh}{Assh}
\DeclareMathOperator{\Min}{Min}
\DeclareMathOperator{\Proj}{Proj}
\DeclareMathOperator{\Ext}{Ext}

\DeclareMathOperator{\height}{height}
\DeclareMathOperator{\codim}{codim}
\DeclareMathOperator{\indeg}{indeg}

\DeclareMathOperator{\depth}{depth}

\DeclareMathOperator{\link}{link}
\newcommand{\frm}{{\mathfrak m}}

\newcommand{\bbZ}{\ensuremath{\mathbb Z}}

\newcommand{\calA}{{\mathcal A}}
\begin{document}
\title{Locally complete intersection Stanley--Reisner ideals}
\author{Naoki Terai}
\address[Naoki Terai]{Department of Mathematics, Faculty of Culture and Education,
         Saga University, Saga 840--8502, Japan}
\email{terai@cc.saga-u.ac.jp}
\thanks{E-mail addresses: \tt terai@cc.saga-u.ac.jp, yoshida@math.nagoya-u.ac.jp}
\author{Ken-ichi Yoshida}
\address[Ken-ichi Yoshida]{Graduate School of Mathematics, Nagoya University,
         Nagoya 464--8602, Japan}
\email{yoshida@math.nagoya-u.ac.jp}
\curraddr{}
\date{\today}
\begin{abstract}
In this paper, we prove that the Stanley--Reisner ideal of 
any connected simplicial complex of dimension $\ge 2$ 
that is locally complete intersection is 
a complete intersection ideal. 
\par 
As an application, we show that the Stanley--Reisner ideal 
whose powers are Buchsbaum is a complete intersection ideal. 
\end{abstract}
\maketitle
\section*{Introduction}
\par 
By a simplicial complex $\Delta$ on a vertex set $V=[n]=\{1,2,\ldots,n\}$,
we mean that $\Delta$ is a non-void family of subsets of $V$
such that (i)$\{v \} \in \Delta $ for every $v \in V$, 
and (ii)$F \in \Delta $, $G \subseteq F$ imply $G \in \Delta $.
Let $S = K[X_1,\ldots,X_n]$ be a polynomial ring over a field $K$. 
The \textit{Stanley--Reisner ideal} of $\Delta$, denoted by $I_{\Delta}$, is 
the ideal of $S$ generated by 
all squarefree monomials $X_{i_1}\cdots X_{i_p}$ such that 
$1 \le i_1 < \cdots < i_p \le n$ and $\{{i_1},\ldots,{i_p}\} \notin \Delta$. 
The \textit{Stanley--Reisner ring} of $\Delta$ over $K$ 
is the $K$-algebra $K[\Delta] = S/I_{\Delta}$.  
Any squarefree monomial ideal $I$ with 
$I \subseteq (X_1,\ldots,X_n)^2$ is a Stanley--Reisner ideal $I_{\Delta}$ 
for some simplicial complex $\Delta$ on $V=[n]$. 
\par
An element $F \in \Delta$ is called a \textit{face} of $\Delta$. 
A maximal face of $\Delta$ with respect to inclusion is called 
a \textit{facet} of $\Delta$. 
The \textit{dimension} of $\Delta$, denoted by $\dim \Delta$, 
is the maximum of the dimensions 
$\dim F = \sharp(F)-1$, where $F$ runs 
through all faces $F$ of $\Delta$ and 
$\sharp(F)$ denotes the cardinality of $F$. 
Note that the Krull dimension of $K[\Delta]$ 
is equal to $\dim \Delta+1$.    
A simplicial complex is called \textit{pure} if all facets 
have the same dimension. 
See \cite{BH, St} for more information on Stanley--Reisner rings.  

\par \vspace{2mm}
A homogeneous ideal $I$ in $S=K[X_1,\ldots,X_n]$ 
is said to be a \textit{locally complete intersection} ideal  
if  $I_P$ is a complete intersection ideal (that is, generated  by a 
regular sequence) for any prime $P \in \Proj (S/I)$. 
A simplicial complex $\Delta$ on $V$ is said to be  
a \textit{locally complete intersection} complex   
if $I_{\link_{\Delta}(\{i\})}$ is a complete intersection ideal 
for every $i \in V$.   
Then $\Delta$ is a locally complete intersection complex
if and only if $I_{\Delta}$ is a locally complete intersection ideal.  
Note that a locally complete intersection ideal $I$ is called  a 
\textit{generalized complete intersection} ideal  
 in the sense of Goto--Takayama (see \cite{GT}) 
if $I=I_{\Delta}$ is the Stanley--Reisner ideal for some 
pure simplicial complex $\Delta$.   

\par \vspace{2mm}
In Section 1, we consider the structure of simplicial complexes which are 
locally complete intersection.  
This is the main purpose of the paper. 
One can easily see that if a Stanley--Reisner ideal $I$ is a 
complete intersection ideal then it 
can be written as 
\[
 I= (X_{11}\cdots X_{1q_1},\ldots, X_{c1} \cdots X_{cq_c}),
\]
where $c \ge 0$ and  $q_i$ is a positive integer with $q_i \ge 2$ for $i=1,\ldots,c$ 
and all $X_{ij}$ are distinct variables.

\par 
A complete intersection simplicial complex $\Delta$ 
is connected if $\dim \Delta \ge 1$, 
and it is a locally complete intersection complex. 
When $\dim \Delta \ge 2$, the converse is also true, 
which is a main result in this paper$:$

\begin{introthm}[\textrm{See also Theorems \ref{Main}, \ref{Structure}}]  
\label{Main-Intro}
Let $\Delta$ be a connected simplicial complex with $\dim \Delta \ge 2$ 
$($resp. $\dim \Delta=1$$)$.
If it is a locally complete intersection complex, 
then it is a complete intersection complex
$($resp. an $n$-gon for $n \ge 3$ or an 
$n$-pointed path for some $n \ge 2$$)$. 
\end{introthm}

\par 
Let $\Delta$ be a connected simplicial complex on $V$ with $\dim \Delta \ge 2$.   
Our main theorem says that 
if $\link_{\Delta} (\{x\})$ is a complete intersection complex 
for every vertex $x \in V$ 
then so is $\Delta$. 
If we also assume Serre's condition $(S_2)$, then we can obtain a stronger result. 
That is, when $K[\Delta]$ satisfies $(S_2)$, 
$\Delta$ is a complete intersection complex if and only if 
$\link_{\Delta}(F)$ is a complete intersection complex 
for any face $F \in \Delta$ with $\dim \link_{\Delta} F =1$;
see Corollary \ref{Link} for more details. 

\par  \vspace{2mm}
In Section 2, we discuss Buchsbaumness for powers 
of Stanley--Reisner ideals. 
Let us explain our motivation briefly.  
Let $A$ be a Cohen--Macaulay local ring. 
If $I$ is a complete intersection ideal of $A$, 
then $A/I^{\ell}$ is Cohen--Macaulay for every $\ell \ge 1$ 
because $I^{\ell}/I^{\ell+1}$ is a free $A/I$-module.  
In \cite{CN}, Cowsik and Nori proved the converse.  
That is, if $I$ is a generically complete intersection ideal
(i.e. $I_P$ is a complete intersection ideal for all minimal 
prime divisors $P$ of $I$) and 
$A/I^{\ell}$ is Cohen--Macaulay for all (sufficiently large) $\ell \ge 1$, 
then $I$ is a complete intersection ideal. 
Note that one can apply this result to Stanley--Reisner ideals:  
$I_{\Delta}$ is a complete intersection ideal if and only if 
$S/I_{\Delta}^{\ell+1}$ is Cohen--Macaulay for every $\ell \ge 1$. 

\par 
A standard graded ring $A=S/I$ with homogeneous 
maximal ideal $\frm$ is said to be \textit{Buchsbaum} (resp. (FLC)) 
if the canonical map 
\[
H^i(\frm,A) \to H_{\frm}^i(A) 
= \varinjlim \Ext_S^i(S/\frm^{\ell},A)
\] 
is surjective (resp. if $H_{\frm}^i(A)$ has finite length) 
for all $i < \dim A$, where 
$H^i(\frm,A)$ (resp. $H^i_{\frm}(A)$) 
denotes the $i$th Koszul cohomology module
(resp. $i$th local cohomology module); 
see \cite[Chapter I, Theorem 2.15]{SV}. 
Then we have the following implications:
\[
\begin{array}{ccccc}
 \text{Complete intersection} 
& \Longrightarrow & \text{Locally complete intersection} & & \\
 & & & & \\
  \Downarrow & & \Downarrow \text{if pure}   & &  \\
 & & & & \\
 \text{Cohen--Macaulay} & \Longrightarrow & \text{Buchsbaum} & 
 \Longrightarrow & \text{(FLC)}. 
\end{array}
\]
\par \vspace{2mm}
Goto and Takayama \cite{GT}  
proved that $I_{\Delta}$ is a pure locally complete intersection ideal 
if and only if 
$S/I_{\Delta}^{\ell+1}$ is (FLC) for every $\ell \ge 1$ 
as an analogue of Cowsik--Nori theorem.

\par
Let $S$ be a polynomial ring and $I$ a squarefree monomial ideal of $S$. 
Then $S/I$ is Buchsbaum if and only if it is (FLC); 
see e.g., \cite[p.73, Theorem 8.1]{St}. 
But a similar statement is no longer true for non-squarefree monomial ideals. 
The following is a natural question$:$

\begin{introquest}  \label{QuestPowers-Intro}
When is $S/I_{\Delta}^{\ell}$ Buchsbaum for every $\ell \ge 1$? 
\end{introquest} 

\par 
As an application of our main theorem and the lower bound formula 
on the multiplicity of Buchsbaum homogeneous $K$-algebras in \cite{GY},  
we can prove the following theorem.  

\begin{introthm}
\label{Powers-Intro}
Put $S = K[X_1,\ldots,X_n]$. 
Let $\Delta$ be a simplicial complex on $V=[n]$. 
Then the following conditions are equivalent$:$
\begin{enumerate}
 \item[(1)] $I_{\Delta}$ is generated by a regular sequence$;$ 
 \item[(2)]  $S/I_{\Delta}^{\ell}$ is Cohen--Macaulay for all $\ell \ge 1$$;$ 
 \item[(3)] $S/I_{\Delta}^{\ell}$ is Buchsbaum for all $\ell \ge 1$$;$  
 \item[(3)'] $\sharp \{\ell \in \bbZ_{\ge 1} \;:\; 
\text{$S/I_{\Delta}^{\ell}$ is Buchsbaum}\}=\infty$.  
\end{enumerate} 
\end{introthm}

\par 
We do not know whether a similar statement is
 true for general homogeneous ideals. 
\bigskip 
\section{Connected complexes which are locally complete intersection}

Throughout this paper, let $\Delta$ be a simplicial complex on $V$. 
For a face $F$ of $\Delta$ and $W \subseteq V$, we put  
\begin{eqnarray*}
\link_{\Delta}(F) &=&
 \{G \in \Delta \;:\; G \cup F \in \Delta,\; F \cap G = \emptyset\}, \\
\Delta_W &=& \{G \in \Delta \;:\; G \subseteq W \}.
\end{eqnarray*}
These complexes are the \textit{link} of $F$, and,  
the \textit{restriction} to $W$ of $\Delta$, respectively. 
\par 
Let $\mathcal{H}$ be a subset of $2^{V}$. 
The minimum simplicial complex $\Gamma \subseteq 2^V$ 
which contains $\mathcal{H}$ as a subset, 
 denoted by $\langle \mathcal{H} \rangle$, 
is said to be the \textit{simplicial complex spanned  
by $\mathcal{H}$ on $V$}.   

\par 
Suppose that $V=V_1 \cup \cdots \cup V_r$ is a disjoint union. 
Let $\Delta_i$ be a simplicial complex on $V_i$ for each $i=1,\ldots,r$. 
Then $\Delta = \Delta_1 \cup \cdots \cup \Delta_r$ is a simplicial complex
on $V$. 
We call $\Delta$ \lq\lq a \textit{disjoint union} of $\Delta_i$'s'' 
by abuse of language although $\Delta_i \cap \Delta_j = \{\emptyset\}$
for $i \ne j$. 

\par \vspace{2mm}
A simplicial complex $\Delta$ is a \textit{complete intersection} complex 
if the Stanley--Reisner ideal $I_{\Delta}$ 
is generated by a regular sequence.  
Now let us define the notion of locally complete intersection 
for complexes.

\begin{defn} \label{LCIdef}
A simplicial complex $\Delta$ on $V$ is 
said to be a \textit{locally complete intersection} complex   
if $I_{\link_{\Delta} (\{i\})}$ 
is a complete intersection ideal  for all vertex $i \in V$. 
\end{defn}

\par 
A simplicial complex $\Delta$ is 
a locally complete intersection complex if and only if  
its Stanley--Reisner ideal $I_{\Delta}$ is 
a locally complete intersection ideal. 

\begin{lemma} \label{Equiv-LCI}
For a Stanley--Reisner ideal $I=I_{\Delta}$, 
the following conditions are equivalent$:$ 
\begin{enumerate}
 \item $\Delta$ is a locally complete intersection complex$;$ 
 \item $K[\Delta]_{X_i}$ is a complete intersection ring for all $i \in V$$;$ 
 \item $I_P$ is a complete intersection ideal for all 
prime $P \in \Proj (S/I_{\Delta})$. 
\end{enumerate}
\end{lemma}
 
\begin{proof}
The equivalence of (1) and (2) immediately follows from the fact that 
\[
 K[\link_{\Delta} (\{i\})][X_i,X_i^{-1}] \cong K[\Delta]_{X_i}.
\] 
(2) $\Longrightarrow (3)$ is clear. 
In order to show the converse, we suppose that 
$K[\Delta]_{X_1}$ is not a complete intersection ring. 
Without loss of generality, we may assume that 
\[
 \{X_i \;:\; 2 \le i \le m\} = \{X_i \;:\; i \in \link_{\Delta}(\{1\}) \}.
\]
Since $X_1X_j \in I_{\Delta}$ for $m+1 \le j \le n$, 
one has that $X_j \in I_{\Delta}S_{X_1}$. 
If we put $P=(X_2,\ldots,X_m)$,  then we can easily see that 
$I_{\Delta}S_P$ is not a complete intersection ideal by assumption.   
Hence we obtain $(3) \Longrightarrow (2)$.  
\end{proof}

\begin{cor} \label{LCIpure}
If $\Delta$ is a connected locally complete intersection complex, 
then it is pure. 
\end{cor}

\begin{proof}
Suppose that $\Delta$ is not pure. 
Since $\Delta$ is connected, there exist a vertex $i \in V$ and 
facets $F_1$, $F_2$ such that $i \in F_1 \cap F_2$ and $\sharp(F_1) < \sharp(F_2)$. 
This implies that $\link_{\Delta}(\{i \})$ is not pure.   
This contradicts the assumption that $\link_{\Delta}(\{i\})$ is Cohen--Macaulay.  
Hence $\Delta$ must be pure. 
\end{proof}

\begin{remark} \label{GCI-comment}
A pure locally 
complete intersection complex is called 
a generalized complete intersection complex in 
\cite{GT}.
\end{remark}

\par 
The main purpose of this section is to prove the following theorem$:$

\begin{thm} \label{Main}
Let $\Delta$ be a connected simplicial complex on $V$ with $\dim \Delta \ge 2$. 
If $\Delta$ is a locally complete intersection complex, 
then it is a complete intersection complex. 
\end{thm}

\par 
Let $\Delta$ be a connected complex of dimension $d-1$.  
Suppose that $\Delta$ is a locally complete intersection complex 
but not a complete intersection complex. 
Note that $\Delta$ is pure and thus a 
generalized complete intersection complex. 
Let $G(I_{\Delta}) = \{m_1,\ldots, m_{\mu}\}$ denote 
the minimal set of monomial generators of $I_{\Delta}$. 
Then $\mu \ge 2$ and $\deg m_i \ge 2$ for every $i =1,2,\ldots,\mu$, and 
that there exist $i$, $j$ ($1 \le i< j \le n$) such that  $\gcd(m_i,m_j) \ne 1$.

\begin{lemma} \label{Twodeg}
In the above notation, we may assume that $\deg m_i  = \deg m_j =2$. 
\end{lemma}

\begin{proof}
Take $m_j$, $m_k$ ($j \ne k$) such that $\gcd(m_j, m_k) \ne 1$. 
If $\deg m_j = \deg m_k=2$, then there is nothing to prove. 
\par 
Now suppose that $\deg m_k \ge 3$. 
By \cite[Lemmas 3.4, 3.5]{GT}, we may assume that $\deg m_j =2$ and 
$\gcd(m_j,m_k) = X_p$. 
Write $m_k = X_p X_{i_1}\cdots X_{i_r}$ and $m_j = X_pX_q$. 
Then \cite[Lemma 3.6]{GT} implies that $X_{i_1}X_q \in G(I_{\Delta})$. 
Set $m_i=X_{i_1}X_q \in I_{\Delta}$. 
Then $\deg m_i = \deg m_j =2$ and $\gcd(m_i,m_j) =X_q \ne 1$, 
as required.  
\end{proof}

\par 
The following lemma is simple but important. 
We use the following convention in this section:
the vertices $x$, $y$, $z$ etc. correspond to the indeterminates $X$, $Y$, $Z$ etc.
respectively. 

\begin{lemma} \label{Three}
Let $x_1,x_2,y$ be distinct vertices such that $X_1Y$, $X_2Y \in I_{\Delta}$.  
For any $z \in V \setminus \{x_1,x_2,y\}$, at lease one of monomials 
$X_1Z$, $X_2Z$ and $YZ$ belongs to $I_{\Delta}$.
\end{lemma}

\begin{proof}
It follows from the fact that $K[\link_{\Delta} (\{z\})]$ is a 
complete intersection ring. 
\end{proof}

\par 
In what follows, we prove Theorem \ref{Main}. 
In order to do that, let $\Delta$ be a connected simplicial complex 
of dimension $d-1 \ge 1$. 
Moreover, assume that $\Delta$ is a locally complete intersection complex 
and that there exist vertices $x_1,x_2,y$ such that 
$X_1Y$, $X_2Y \in I_{\Delta}$
(we assign a variable $X_i$ for a vertex $x_i$).
Then we must show that $\dim \Delta (=d-1) =1$.  
Let us begin with proving the following key lemma. 

\begin{lemma} \label{Key}
Under the above notation, 
there exist some integers $k$, $\ell \ge 2$ such that 
\begin{enumerate}
 \item $V = \{x_1,\ldots,x_k,y_1,\ldots,y_{\ell} \}$$;$ 
 \item $X_1Y_1,\ldots,X_kY_1 \in I_{\Delta}$$;$ 
 \item $\sharp\{i \;:\; 1 \le i \le k,\; X_iY_j \notin I_{\Delta} \} \le 1$ holds 
for each $j=2,\ldots, \ell$. 
\end{enumerate}
\end{lemma}

\begin{proof}
By assumption, there exist vertices $x_1,x_2,y_1 \in V$ 
such that $X_1Y_1$, $X_2Y_1 \in I_{\Delta}$. 
Thus one can write $V = \{x_1,\ldots,x_k,y_1,\ldots,y_{\ell}\}$ such that 
\[
\begin{array}{cccccc}
X_1Y_1, & X_2Y_1, & \ldots, & X_kY_1 & \in & I_{\Delta}, \\
Y_1Y_2, & Y_1Y_3, & \ldots, & Y_1Y_{\ell} & \notin & I_{\Delta}. \\
\end{array}
\]
If $\ell =1$, then 
$\Delta = \Delta_{\{y_1\}} \cup \Delta_{\{x_1,\ldots,x_k\}}$ 
is a disjoint union 
since $\{y_1,x_i\} \notin \Delta$ for all $i$. 
This contradicts the connectedness of $\Delta$. 
Hence $\ell \ge 2$. 
Thus it is enough to show (3) in this notation. 

\par 
Now suppose that there exists an integer $j$ with $2 \le j \le \ell$ 
such that 
\[
 \sharp\{i \;:\; 1 \le i \le k,\, X_iY_j \notin I_{\Delta} \} \ge 2. 
\]
When $k=2$, we have $X_1Y_j$, $X_2Y_j \notin I_{\Delta}$. 
On the other hand, as $X_1Y_1$, $X_2Y_1 \in I_{\Delta}$ and $Y_j \ne X_1,X_2,Y_1$, 
we obtain that at least one of $X_1Y_j$, $X_2Y_j$, $Y_1Y_j$ belongs to $I_{\Delta}$. 
It is impossible. 
So we may assume that $k \ge 3$ and $X_{k-1}Y_j$, $X_k Y_j \notin I_{\Delta}$. 
Then $\{x_{k-1}\}$, $\{x_k\}$ and $\{y_1\}$ belong to $\link_{\Delta}(\{y_j\})$, and 
$X_{k-1}Y_1$, $X_kY_1$ form part of the minimal system of 
generators of $I_{\link_{\Delta} (\{y_j\})}$. 
This contradicts the assumption that $\link_{\Delta} (\{y_j\})$ is a 
complete intersection complex. 
\end{proof}

\par
In what follows, we fix the notation as in Lemma  \ref{Key}. 
First, we suppose that there exists an $i_0$ with $1 \le i_0 \le k$ such that 
\[
 \sharp\{j \;:\;1 \le j \le \ell,\, X_{i_0}Y_j \notin I_{\Delta} \}=1.
\] 
In this case, we may assume that $X_1Y_2 \notin I_{\Delta}$ and 
$X_1Y_j \in I_{\Delta}$ for all $3 \le j \le \ell$ without loss of generality. 
Note that $X_2Y_2,\ldots, X_{k}Y_2 \in I_{\Delta}$ by Lemma \ref{Key}. 
We claim that $\{x_1,y_2\}$ is a facet of $\Delta$. 
As $X_{i}Y_2 \in I_{\Delta}$ for each $i=2,\ldots,k$, 
it follows that $\{x_1,y_2,x_i\} \notin \Delta$. 
Similarly, $\{x_1,y_2,y_j\} \notin \Delta$ since 
$X_1Y_j \in I_{\Delta}$ for $j=1$ or $3 \le j \le \ell$. 
Hence $\{x_1,y_2\}$ is a facet of $\Delta$, and 
$\dim \Delta =1$ because $\Delta$ is pure. 

\par\vspace{2mm}
By the observation as above, we may assume that for every $i$ with $1 \le i \le k$,
\[
  \sharp\{j \;:\;1 \le j \le \ell,\, X_{i}Y_j \notin I_{\Delta} \} \ge 2
\]
or $X_iY_j \in I_{\Delta}$ holds for all $j=1,\ldots,\ell$. 
\par
Now suppose that there exist $j_1,j_2$ with $1 \le j_1 < j_2 \le \ell$ such that 
$X_iY_{j_1}$, $X_{i}Y_{j_2} \notin I_{\Delta}$. 
Then $X_{r}Y_{j_1}$, $X_rY_{j_2} \in I_{\Delta}$ for all $r \ne i$ by Lemma \ref{Key}. 
It follows that $X_rX_i \in I_{\Delta}$ from Lemma \ref{Three}. 
Then we can relabel $x_i$ (say $y_{\ell+1})$. 
Repeating this procedure, we can get one of the following cases:

\begin{description}
\item[\textbf{Case 1}] $V=\{x_1,\ldots,x_r,y_1,\ldots,y_s\}$ such that 
$X_iY_j \in I_{\Delta}$ for all $i$, $j$ 
with $1 \le i \le r$, $1 \le j \le s$. 
\vspace{2mm}
\item[\textbf{Case 2}] 
$V = \{x_1,x_2,y_1,\ldots,y_m,z_1,\ldots,z_p,w_1,\ldots,w_q\}$ 
such that 
\[
\left\{
\begin{array}{ccc}
X_1Y_j \in I_{\Delta}, & X_2Y_j \in I_{\Delta} & (j=1,\ldots,m)\\
X_1Z_j \notin I_{\Delta}, & X_2Z_j \in I_{\Delta} & (j=1,\ldots,p)\\
X_1W_j \in I_{\Delta}, & X_2W_j \notin I_{\Delta} & (j=1,\ldots,q)
\end{array}
\right.
\]
holds for some $m \ge 1$, $p,q \ge 2$. 
\end{description}
If Case 1 occurs, then 
$\Delta = \Delta_{\{x_1,\ldots,x_r\}} \cup \Delta_{\{y_1,\ldots,y_s\}}$ is 
a disjoint union. 
This contradicts the assumption. 
Thus Case 2 must occur. 
If $\{x_1,x_2\} \in \Delta$, then it is a facet and so $\dim \Delta=1$. 
Hence we may assume that $\{x_1,x_2\} \notin \Delta$.  
However, since $\Delta$ is connected, there exists a path between $x_1$ and $x_2$. 

\vspace{2mm}
\begin{description}
\item[{\bf Cases (2-a)}] the case where  $\{z_1,w_k\} \in \Delta$ for some $k$ with $1 \le k \le q$. 
\end{description}
\par 
We may assume that $\{z_1,w_1\} \in \Delta$.  
Now suppose that $\dim \Delta \ge 2$. 
Then since $\{z_1,w_1\}$ is \textit{not} a facet,  
there exists a vertex $u \in V \setminus \{x_1,x_2\}$ 
such that $\{z_1,w_1,u \} \in \Delta$. 
If $u=z_j$ $(2 \le j \le p)$ (resp. $u=y_i$ $(1 \le i \le m)$), then   
$G(I_{\link_{\Delta} (\{w_1\})})$ contains $X_2Z_1$ and $X_2Z_j$
(resp. $X_2Y_i$); see the figure below.  
It is impossible since $\link_{\Delta} (\{w_1\})$ is a complete intersection complex.   
When $u=w_k$, we can obtain a contradiction by a similar argument as above. 
Therefore $\dim \Delta=1$. 
\vspace{4mm}
\begin{center}
\begin{picture}(160,40)
  \put(0,15){\circle*{5}}
  \put(40,34){\circle*{5}}
  \put(40,-4){\circle*{5}}
  \put(80,34){\circle*{5}}
  \put(120,15){\circle*{5}}
  \put(0,15){\line(2,1){37}}
  \put(0,15){\line(2,-1){37}}
  \put(120,15){\line(-2,1){37}}
  \put(40,34){\line(1,0){40}}
  \put(40,34){\line(0,-1){38}}
  \put(40,-4){\line(1,1){40}}
  \put(-15,14){$x_1$}
  \put(38,41){$z_1$}
  \put(38,-13){$z_j$}
  \put(80,41){$w_1$}
  \put(126,14){$x_2$}
  \put(52,21){\circle{8}}
\end{picture}
\par \vspace{6mm}
Figure: the case $\{z_1,z_j,w_1\} \in \Delta$ in Case (2-a)
\end{center}
\vspace{5mm}
\begin{description}
\item[{\bf Cases (2-b)}] the case where  $\{z_j,w_k\} \notin \Delta$ 
for all $j$, $k$. 
\end{description}
\par 
Then we may assume that (i) $\{z_1,y_1\} \in \Delta$ 
and (ii) $\{y_1,y_2\} \in \Delta$ or $\{y_1,w_1\} \in \Delta$. 
Now suppose that $\dim \Delta \ge 2$. 
Then since $\{z_1,y_1\}$ is \textit{not} a facet,  we have 
\[
\{z_1, y_1, y_i\} \in \Delta, \;   
\{z_1, y_1,w_k\} \in \Delta \; \text{or}\;
\{z_1, y_1, z_j\} \in \Delta. 
\]
When $\{z_1,y_1, y_i \} \in \Delta$, we obtain that 
$X_1Y_1,X_1Y_i \in G(I_{\link_{\Delta} (\{z_1\})})$. 
This is a contradiction. 
When $\{z_1,y_1,w_k\} \in \Delta$, 
we can obtain a contradiction by a similar argument as in Case (2-a). 
Thus it is enough to consider the case $\{z_1,y_1,z_j \} \in \Delta$.  

\par \vspace{2mm}
First we suppose that $\{y_1,y_2\} \in \Delta$. 
\vspace{4mm}
\begin{center}
\begin{picture}(160,40)
  \put(0,15){\circle*{5}}
  \put(40,34){\circle*{5}}
  \put(40,-4){\circle*{5}}
  \put(80,34){\circle*{5}}
  \put(130,15){\circle*{5}}
  \put(120,34){\circle*{5}}
  \put(0,15){\line(2,1){37}}
  \put(0,15){\line(2,-1){37}}
  \put(80,34){\line(1,0){40}}
  \put(40,34){\line(1,0){40}}
  \put(40,34){\line(0,-1){38}}
  \put(40,-4){\line(1,1){40}}
  \put(-15,14){$x_1$}
  \put(38,41){$z_1$}
  \put(38,-13){$z_j$}
  \put(80,41){$y_1$}
  \put(118,41){$y_2$}
  \put(136,14){$x_2$}
  \put(52,21){\circle{8}}
\end{picture}
\par \vspace{6mm}
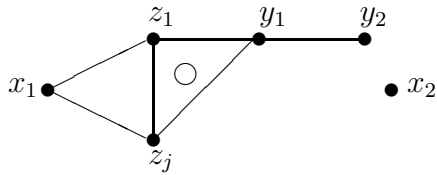
Figure: the case $\{z_1,y_1,z_j\}, \{y_1,y_2\} \in \Delta$ in Case (2-b)
\end{center}
\vspace{5mm}
Then $\link_{\Delta} (\{y_1\})$ contains  $\{z_1,z_j\}$ and $\{y_2\}$. 
Since $\link_{\Delta} (\{y_1\})$ is also connected, 
we can find vertices $z_\alpha$, $y_\beta$ such that 
$\{z_{\alpha},y_{\beta}\} \in \link_{\Delta} (\{y_1\})$. 
In particular, $\{z_\alpha,y_{\beta},y_1\} \in \Delta$. 
This yields a contradiction because 
$X_1Y_1, X_1Y_{\beta}$ are contained in $G(I_{\link_{\Delta}(\{z_{\alpha}\})})$. 
\par \vspace{2mm}
Next suppose that $\{y_1,w_1\} \in \Delta$.
\vspace{4mm}
\begin{center}
\begin{picture}(160,40)
  \put(0,15){\circle*{5}}
  \put(40,34){\circle*{5}}
  \put(40,-4){\circle*{5}}
  \put(80,34){\circle*{5}}
  \put(136,18){\circle*{5}}
  \put(120,34){\circle*{5}}
  \put(0,15){\line(2,1){37}}
  \put(0,15){\line(2,-1){37}}
  \put(80,34){\line(1,0){40}}
  \put(40,34){\line(1,0){40}}
  \put(40,34){\line(0,-1){38}}
  \put(40,-4){\line(1,1){40}}
  \put(120,34){\line(1,-1){15}}
  \put(-15,14){$x_1$}
  \put(38,41){$z_1$}
  \put(38,-13){$z_j$}
  \put(80,41){$y_1$}
  \put(118,41){$w_1$}
  \put(140,16){$x_2$}
  \put(52,21){\circle{8}}
\end{picture}
\par \vspace{6mm}
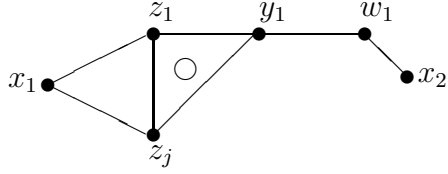
Figure: the case $\{z_1,y_1,z_j\}, \{y_1,w_1\} \in \Delta$ in Case (2-b)
\end{center}
\vspace{5mm}
Then $\link_{\Delta} (\{y_1\})$ contains $\{z_1,z_j\}$ and $\{w_1\}$. 
Since $\link_{\Delta} (\{y_1\})$ is also connected, 
we can also find vertices 
$z_\alpha$, $y_\beta$ such that 
$\{z_{\alpha},y_{\beta}\} \in \link_{\Delta} (\{y_1\})$
(notice that $\{z_j,w_k\} \notin \Delta$). 
In particular, $\{z_{\alpha},y_1,y_{\beta}\} \in \Delta$. 
This yields a contradiction because 
$X_1Y_1, X_1Y_{\beta}$ are contained in $G(I_{\link_{\Delta}(\{z_{\alpha}\})})$. 

\par \vspace{2mm}
Therefore we have $\dim \Delta=1$. 
So we have finished the proof of Theorem \ref{Main}. 
 
\par \vspace{3mm}
An arbitrary Noetherian ring $R$ is said to satisfy
Serre's condition $(S_2)$ if  
$\depth R_P \ge \min\{\dim R_P,\, 2\}$ for every prime $P$ of $R$. 
A Stanley--Reisner ring $K[\Delta]$ satisfies $(S_2)$ if and only if  
$\Delta$ is pure and $\link_{\Delta}(F)$ is connected for every face $F$ with 
$\dim \link_{\Delta}(F) \ge 1$; see e.g., \cite[p.454]{Te}. 
In particular, if $K[\Delta]$ satisfies $(S_2)$, then $\Delta$ is 
pure and connected if $\dim \Delta \ge 1$. 

\par 
Let $\Delta$ be a connected simplicial complex on $V$ with $\dim \Delta \ge 2$.   
Our main theorem says that 
if $\link_{\Delta} (\{x\})$ is a complete intersection complex for every $x \in V$ 
then so is $\Delta$ itself. 
Thus it is natural to ask the following question:

\begin{quest}  \label{Subset-Intro}
Does there exist a proper subset $W \subseteq V$ for which  
\lq\lq \, $\link_{\Delta} (\{x\})$ is 
a complete intersection complex for all $x \in W$'' 
implies that $\Delta$ is a complete intersection complex?
\end{quest}

\par 
The following corollary gives an answer to the above question 
in the $(S_2)$ case.

\begin{cor} \label{Link}
Let $\Delta$ be a simplicial complex with $\dim \Delta \ge 2$. 
Assume that $K[\Delta]$ satisfies $(S_2)$. 
Then the following conditions are equivalent$:$
\begin{enumerate}
 \item[(1)] $K[\Delta]$ is a complete intersection ring$;$ 
 \item[(2)] For any face $F$ with $\dim \link_{\Delta} (F) =1$, 
$\link_{\Delta} (F)$ is a complete intersection complex$;$  
 \item[(3)] There exists $W \subseteq V$ such that 
$\dim \Delta_{V \setminus W} \le \dim \Delta -3$ 
which satisfies the following condition$:$
\par  \vspace{1mm}
\lq\lq \, $\link_{\Delta} (\{x\})$ is a complete intersection complex
 for all $x \in W$.''. 
\end{enumerate} 
\end{cor}

\begin{proof}
Note that $\Delta$ is pure. 
Put $d = \dim \Delta+1$. 
\par \vspace{2mm}  
$(1) \Longrightarrow (3):$ It is enough to put $W=V$. 
\par \vspace{2mm}  
$(3) \Longrightarrow (2):$ Let $W \subseteq V$ be a subset that satisfies 
the condition (3). 
Let $F$ be a face with $\dim \link_{\Delta}(F) =1$. 
Since $\Delta$ is pure, $\sharp(F) = d-1-\dim \link_{\Delta} (F) = d-2$.
As $\dim \Delta_{V \setminus W}  \le d-4$, $F$ is not contained in $V \setminus W$. 
Thus there exists a vertex $i \in F$ such that $i \in W$. 
Then since $\link_{\Delta} (\{i\})$ is a complete intersection complex 
by assumption, 
$\link_{\Delta}(F)$ is also a complete intersection complex, as required.  
\par \vspace{2mm}  
$(2) \Longrightarrow (1):$ 
We use an induction on $d \ge 3$. 
First suppose that $d=3$. 
Then for each $i \in V$, since $\dim \link_{\Delta} (\{i\}) =1$, 
$\link_{\Delta}(\{i\})$ is a complete intersection complex by the assumption (2). 
Hence $K[\Delta]$ is a complete intersection ring by Theorem \ref{Main}. 
\par 
Next suppose that $d \ge 4$. 
Let $i \in V$. 
Since $K[\Delta]$ satisfies $(S_2)$, we have that 
$\Gamma=\link_{\Delta} (\{i\})$ is connected and 
$\dim \Gamma =d-2 \ge 2$. 
Moreover, for any face $G$ in $\Gamma$ with $\dim \link_{\Gamma} (G)=1$, 
$\link_{\Gamma} (G)=\link_{\Delta}(G \cup \{i\})$ is a complete intersection complex 
by assumption.   
Hence, by the induction hypothesis, 
$K[\link_{\Delta}(\{i\})]$ is a complete intersection ring. 
Therefore $K[\Delta]$ is a complete intersection ring 
by Theorem \ref{Main} again.
\end{proof} 

\par 
To complete the proof of Theorem \ref{Main-Intro}, we must consider the case 
$\dim \Delta=1$. 
In this case, there exist connected non-complete intersection complexes 
that are locally complete intersection.
\par 
Let $\Delta$ be a one-dimensional simplicial complex on $V=[n]$. 
$\Delta$ is said to be the \textit{$n$-gon} for $n \ge 3$ 
(resp. the \textit{$n$-pointed path} for $n \ge 2$) if 
$\Delta$ is pure and its facets consist of $\{i,i+1\}$ 
$(i=1,2\ldots,n-1)$ and $\{n,1\}$ 
(resp.  its facets consists of $\{i,i+1\}$ $(i=1,2\ldots,n-1)$)
after suitable change of variables. 

\begin{prop} \label{1dimLCI} 
Let $\Delta$ be a $1$-dimensional connected complex.  
Then the following conditions are equivalent$:$ 
\begin{enumerate}
 \item $\Delta$ is a locally complete intersection complex$;$ 
 \item $\Delta$ is locally Gorenstein $($i.e., $K[\link_{\Delta}(\{i\})]$ 
is Gorenstein for every $i \in V$$);$ 
 \item $\Delta$ is isomorphic to either one of the following$:$
\begin{enumerate}
 \item the $n$-gon for $n \ge 3$$;$
 \item the $n$-pointed path for $n \ge 2$. 
\end{enumerate}
\end{enumerate}
\end{prop}

\begin{proof}
Note that $(1) \Longrightarrow (2)$ is clear. 
\par 
Suppose that $\Delta$ is a locally Gorenstein. 
Then since $\link_{\Delta}(\{i\})$ is a zero-dimensional 
Gorenstein complex, it consists of at most two points.  
Such a complex is isomorphic to either one of the $n$-gon ($n \ge 3$) or 
the $n$-pointed path ($n \ge 2$). 
\par 
Conversely, if $\Delta$ is isomorphic to either $n$-gon or $n$-pointed path, 
then $\link_{\Delta} (\{i\})$ is a complete intersection complex. 
Hence $\Delta$ is locally complete intersection. 
\end{proof}

\begin{remark} \label{LCIbutNotCI}
Let $\Delta$ be a connected simplicial complex 
on $V=[n]$ of $\dim \Delta=1$. 
Then $\Delta$ is a locally complete intersection complex
but not a complete intersection complex 
if and only if  it is isomorphic to the $n$-gon 
for some $n \ge 5$ or the $n$-pointed path 
for some $n \ge 4$. 
\end{remark}

\begin{exam}  \label{Gon}
Let $K$ be a field. 
The Stanley--Reisner ring of the $4$-pointed path $\Delta_1$ 
is 
$K[X_1,X_2,X_3,X_4]/(X_1X_3,X_1X_4,X_2X_4)$. 
The Stanley--Reisner ring of the $5$-gon $\Delta_2$ is 
$K[X_1,X_2,X_3,X_4,X_5]/(X_1X_3,X_1X_4,X_2X_4,X_2X_5,X_3X_5)$.

\par \noindent
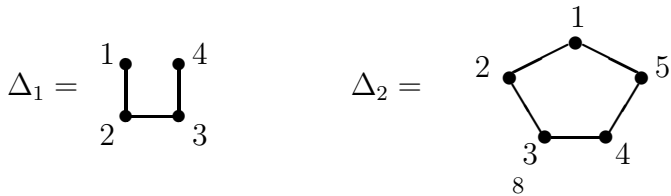
\begin{picture}(400,70)
 \put(20,25){$\Delta_1=$}
  \thicklines
\put(65,17){\line(1,0){20}} 
\put(65,15){\line(0,1){20}} 
\put(85,15){\line(0,1){20}} 
\put(62,14){$\bullet$}
\put(82,14){$\bullet$}
\put(62,34){$\bullet$}
\put(82,34){$\bullet$}
\put(55,37){$1$}
\put(55,6){$2$}
\put(90,6){$3$}
\put(90,37){$4$}
\put(150, 25){$\Delta_2=$}
  \thicklines
  \put(235,45){\circle*{5}}  
  \put(233,51){$1$}
  \put(210,31.5){\circle*{5}}  
  \put(197,32){$2$}
  \put(223.5,9){\circle*{5}}  
  \put(215,-1){$3$}
  \put(246.5,9){\circle*{5}}  
  \put(250,-1){$4$}
  \put(260,31.5){\circle*{5}}  
  \put(265,32){$5$}
  \put(232,44){\line(-2,-1){20.5}}  
  \put(238,44){\line(2,-1){20.5}}  
  \put(222,11){\line(-3,5){10.8}}  
  \put(248,11){\line(3,5){10.8}}  
  \put(226,9){\line(1,0){18}}  
\end{picture}
\end{exam}

\begin{remark}
When $\dim \Delta \ge 2$, 
there are many examples of locally Gorenstein complexes 
which are not locally complete intersection complexes. 
\end{remark}  

\par 
In the last of this section, we give a structure theorem for locally 
complete intersection complexes. 

\begin{thm} \label{Structure}
Let $\Delta$ be a simplicial complex on $V$ such that $V \ne \emptyset$. 
Then $\Delta$ is a locally complete intersection complex if and only if 
it is a finitely many disjoint union of the following connected 
complexes$:$ 
\begin{enumerate}
\item[(a)] a complete intersection complex $\Gamma$ with $\dim \Gamma \ge 2;$ 
\item[(b)] $m$-gon $(m \ge 3);$
\item[(c)] $m'$-pointed path $(m' \ge 2)$. 
\item[(d)] a point 
\end{enumerate}
When this is the case, $K[\Delta]$ is Cohen--Macaulay 
$($resp. Buchsbaum $)$ 
if and only if $\dim \Delta =0$ or $\Delta$ is connected  $($resp. pure$)$. 
\end{thm} 

\par 
To prove the theorem, it suffices to show the following lemma. 

\begin{lemma} \label{Discon}
Assume that $V = V_1 \cup V_2$ such that $V_1 \cap V_2 = \emptyset$. 
Let $\Delta_i$ be a simplicial complex on $V_i$ for $i=1,2$.  
If $\Delta_1$ and $\Delta_2$ are both locally complete intersection complexes, 
then so is $\Delta_1 \cup \Delta_2$. 
\end{lemma}

\begin{proof}
Put $\Delta = \Delta_1 \cup \Delta_2$ and 
$V_1 = [m]$ and $V_2 = [n]$. 
If we write 
\[
K[\Delta_1] =K[X_1,\ldots,X_m]/I_{\Delta_1}
\quad 
\text{and}
\quad 
K[\Delta_2] =K[Y_1,\ldots,Y_n]/I_{\Delta_2}, 
\]
then 
\[
 K[\Delta] \cong K[X_1,\ldots,X_m,Y_1,\ldots,Y_n]/
(I_{\Delta_1},\, I_{\Delta_2},\,\{X_iY_j\}_{1\le i \le m,\; 1\le j \le n}).
\]
Hence $K[\Delta]_{X_i} \cong K[\Delta_1]_{X_i}$ and  
$K[\Delta]_{Y_j} \cong K[\Delta_2]_{Y_j}$
are complete intersection rings. 
Thus $\Delta$ is also a locally complete intersection 
complex by Lemma \ref{Equiv-LCI}. 
\end{proof}

\begin{remark} \label{Gen}
In the above lemma, we suppose that 
both $\Delta_1$ and $\Delta_2$ are 
generalized complete intersection complexes. 
Then $\Delta_1 \cup \Delta_2$ is a generalized complete intersection complexes
if and only if $\dim \Delta_1 = \dim \Delta_2$.   
\end{remark}

\begin{exam} \label{Ex-disjoint}
Let $\Delta$ be the disjoint union of the standard $(m-1)$-simplex 
and the standard $(n-1)$-simplex. 
Then $\Delta$ is a locally complete intersection complex 
by Lemma \ref{Discon}. 
Moreover, $K[\Delta]$ is isomorphic to  
\[
 K[X_1,\ldots,X_m,Y_1,\ldots,Y_n]/(X_iY_j \;:\; 1 \le i \le m,\, 1 \le j \le n) 
\]
and it is a generalized complete intersection complex 
if and only if $m=n$.  
\end{exam}


\medskip
\section{Buchsbaumness of powers for Stanley--Reisner ideals}

\par 
The Stanley--Reisner ring $K[\Delta]$ has (FLC) 
if and only if $\Delta$ is pure and 
$K[\link_{\Delta} (\{i\})]$ is Cohen--Macaulay for every $i \in V$. 
Then $H_{\frm}^i(K[\Delta]) = [H_{\frm}^i(K[\Delta])]_0$ for all $i < \dim K[\Delta]$ 
and  so that $K[\Delta]$ is Buchsbaum. 
See \cite[p.73, Theorem 8.1]{St}.   

\par 
Let $\ell \ge 2$ be an integer. 
Suppose that $S/I_{\Delta}^{\ell}$ is Buchsbaum. 
In \cite{HTT}, Herzog, Takayama and the first author showed that 
this condition implies that $S/I_{\Delta}$ is Buchsbaum.   
The converse is not true. 
What can we say about the structure of $\Delta$?  
This gives a motivation of our study in this section. 

\par \vspace{2mm}
The main result in this section is the following theorem, which is an analogue
of the Cowsik--Nori theorem in \cite {CN} and 
the Goto--Takayama theorem in \cite{GT}.  

\begin{thm} \label{BbmPowers} 
Put $S = K[X_1,\ldots,X_n]$. 
Let $I_{\Delta}$ denote the Stanley--Reisner ideal 
of a simplicial complex $\Delta$ on $V=[n]$. 
Then the following conditions are equivalent$:$
\begin{enumerate}
 \item[(1)] $I_{\Delta}$ is generated by a regular sequence$;$ 
 \item[(2)]  $S/I_{\Delta}^{\ell}$ is Cohen--Macaulay for all $\ell \ge 1$$;$ 
 \item[(3)] $S/I_{\Delta}^{\ell}$ is Buchsbaum for all $\ell \ge 1$$;$ 
 \item[(3)'] $\sharp \{\ell \in \bbZ_{\ge 1} \;:\; 
\text{$S/I_{\Delta}^{\ell}$ is Buchsbaum}\}=\infty$.  
\end{enumerate} 
\end{thm}

\par
Note that $(1) \Longleftrightarrow (2)$ is a special case of 
the Cowsik--Nori theorem and 
$(2) \Longrightarrow (3) \Longrightarrow (3)'$ is trivial. 
Thus our contribution is $(3)' \Longrightarrow (1)$.  

\par \vspace{3mm}
In what follows, 
we put $d = \dim S/I_{\Delta}$, 
$c=\height I_{\Delta} (= \codim I_{\Delta})=n-d$.  
Put $q=\indeg I_{\Delta} \ge 2$, the \textit{initial degree} of $I$, 
that is, $q$ is the least degree of the minimal generators of $I$,
in other words, $q=\min\{\sharp(F) \;:\; F \in 2^{V} \setminus \Delta \}$.
Put $e=e(S/I_{\Delta})$, the \textit{multiplicity} of $I_{\Delta}$, 
which is equal to the number of facets of dimension $d-1$. 
Note that for any homogeneous ideal $I$ of $S$, the following 
formula for multiplicities is known$:$
\[
 e(S/I) = \sum_{P \in \Assh_S(S/I)} e(S/P) \cdot \lambda_{S_P}(S_P/IS_P),
\] 
where $\Assh_S(S/I) = \{P \in \Min_S(S/I) \,:\, \dim S/P = \dim S/I\}$
and $\lambda_R(M)$ denotes the length of an $R$-module $M$ over 
an Artinian local ring $R$. 
\par   
In order to prove the theorem, it suffices to show that  
if $S/I_{\Delta}^{\ell}$ is Buchsbaum for infinitely many $\ell \ge 1$, 
then  $\Delta$ is a complete intersection complex. 
\par 
First we give a formula for multiplicities of  $S/I_{\Delta}^{\ell}$ 
for every $\ell \ge 1$. 

\begin{lemma} \label{formula}
Under the above notation, we have 
\[
 e(S/I_{\Delta}^{\ell}) = e \cdot \binom{c+\ell-1}{c}. 
\]
\end{lemma}

\begin{proof}
Let $P \in \Assh_S(S/I_{\Delta}^{\ell})$. 
Then $P$ is a minimal prime over $I_{\Delta}$ such that 
$S/P$ is isomorphic to a polynomial ring in $d$ variables and 
$S_P$ is a regular local ring of dimension $c$. 
Thus we get
\[ 
 e(S/I_{\Delta}^{\ell}) = \sum_{P \in \Assh S/I_{\Delta}} \!\! e(S/P) 
\cdot \lambda_{S_P} (S_P/I_{\Delta}^{\ell}S_P) 
= e \cdot \binom{c+\ell-1}{c}, 
\]
as required. 
\end{proof}

\par 
We recall the following theorem, which gives a lower bound on multiplicities 
for homogeneous Buchsbaum algebras$:$

\begin{lemma}[{\rm \cite[Theorem 3.2]{GY}}]
Assume that $S/I$ is a homogeneous Buchsbaum $K$-algebra. 
Put $c=\codim I \ge 2$, $q=\indeg I \ge 2$ and $d = \dim S/I \ge 1$. 
Then 
\[
 e(S/I) \ge \binom{c+q-2}{c} + \sum_{i=1}^{d-1} 
\binom{d-1}{i-1} \cdot \dim_K H_{\frm}^i(S/I). 
\]  
\end{lemma}

\par \vspace{2mm}
Applying this formula to $S/I^{\ell}_{\Delta}$ yields 

\begin{cor} \label{Mult-Bbm}
If $S/I_{\Delta}^{\ell}$ is Buchsbaum, then 
\[
 e(S/I_{\Delta}^{\ell}) \ge \binom{c+q\ell-2}{c}.  
\]  
In particular, we have 
\[
 e(S/I_{\Delta}) \ge 
\dfrac{\binom{c+q\ell-2}{c}}{\binom{c+\ell-1}{c}}  
= \dfrac{(q\ell+c-2)\cdots (q\ell+1)q\ell(q\ell-1)}{(\ell+c-1)\cdots (\ell+1)\ell}. 
\]
\end{cor}

\par \vspace{2mm}
In the above corollary, 
if we fix $c, q$ and let $\ell$ tend to $\infty$, then the limit of the right hand side 
in the last inequality tends to $q^c$.  
Therefore if $S/I_{\Delta}^{\ell}$ is Buchsbaum for infinitely many $\ell \ge 1$, 
then $e(S/I_{\Delta}) \ge q^c$.  
For instance, if $I_{\Delta}=(m_1,\ldots,m_{c})$ is a complete intersection ideal, 
then this inequality holds because 
\[
 e(S/I_{\Delta}) = \deg m_1 \cdots \deg m_c \ge q^c.  
\]
However, if $I$ is a locally complete intersection ideal 
but not a complete intersection ideal, 
then this is not true. 
This is a key point in the proof of Theorem \ref{BbmPowers}. 
Namely we have$:$ 

\begin{prop} \label{LCI-mult}
Assume that $\Delta$ is pure and a locally complete intersection complex  
but not a complete intersection complex. 
Then
\[ 
 e(K[\Delta]) < 2^c. 
\]   
\end{prop}

\begin{proof}
First we consider the case $d=1$. 
Then $\Delta$ consists of $n$ points, and so that $c=n-1$, $e=n$. 
As $\Delta$ is not a complete intersection complex, we have $n \ge 3$. 
Then $e=n < 2^c= 2^{n-1}$ is clear. 
\par 
Next we consider the case $d=2$. 
By assumption, $\Delta$ is isomorphic to the following complexes:
\begin{enumerate}
 \item[(a)] the $n$-gon for $n \ge 5$;   
 \item[(b)] the $n$-pointed path for $n \ge 4$; 
 \item[(c)] the disjoint union of $k$ connected 
complexes $\Delta_1, \ldots, \Delta_k$ for some $k \ge 2$, where 
each $\Delta_i$ is isomorphic to the $m$-gon for some $m \ge 3$ 
or the $m$-pointed path for $m \ge 2$. 
\end{enumerate} 
\par 
In particular, we have $e \le n$ and $c=n-2$. 
If $n \ge 5$, then $e \le n < 2^{n-2}=2^c$ is clear.   
So we may assume that $3 \le n \le 4$. 
Then $\Delta$ is isomorphic to either the $4$-pointed path or 
two disjoint union of the $2$-pointed paths. 
In any case, we have $e \le 3 < 4=2^{c}$. 
\par 
Finally, we consider the case $d \ge 3$. 
Theorem \ref{Main} implies that $\Delta$ is disconnected, and so that $c \ge d$. 
Then we consider the following three cases: 
\begin{enumerate}
 \item[(a)] the case $c=d;$
 \item[(b)] the case $c=d+1;$
 \item[(c)] the case $c \ge d+2$. 
\end{enumerate}
\par 
When $c=d$, $\Delta$ is a disjoint union of two $(d-1)$-simplices.  
Then $e = 2 < 2^3 \le 2^c$, as required. 
When $c=d+1$, $\Delta$ has just two connected components. 
One of components is a $(d-1)$-simplex and 
the other one is a pure $(d-1)$-subcomplex of the boundary complex of a $d$-simplex. 
In particular, $e \le d+2 < 2^{c}=2^{d+1}$.  
\par 
So we may assume that $c \ge d+2$. 
Then $\Delta$ is a disjoint union of complete intersection complexes 
of dimension $d-1$ (say, $\Delta_1$,\ldots, $\Delta_k$) by Theorem \ref{Structure}, 
where $k \le \frac{n}{d}=1+ \frac{c}{d}$. 
Moreover, since $c \ge d+2$, we obtain that 
$c(d-1) \ge (d+2)(d-1) > d^2$, and thus $d+\frac{c}{d} <c$. 
Hence 
\[
 e(K[\Delta]) = \sum_{i=1}^k e(K[\Delta_i]) 
\le 2^d \cdot k \le 2^d \cdot (1+\frac{c}{d}) \le 2^d \cdot 2^{\frac{c}{d}} 
=2^{d+\frac{c}{d}} < 2^c, 
\]  
where the first inequality follows from the lemma below. 
\end{proof} 

\begin{lemma} \label{CI-mult}
Assume that $\Delta$ is a complete intersection complex of dimension $d-1$. 
Then $e(K[\Delta]) \le 2^d$.  
\end{lemma}

\begin{proof}
Write $I_{\Delta} = (m_1,\ldots,m_{c})$, where $\deg m_i = h_i$ $(i=1,\ldots,c)$. 
Then 
\[
 e(K[\Delta]) = h_1\cdots h_c \le 2^{h_1-1} \cdots 2^{h_c-1} 
= 2^{h_1+\cdots + h_c-c} \le 2^{n-c} = 2^d, 
\]
as required. 
\end{proof}

\par 
We are now ready to prove Theorem \ref{BbmPowers}. 
 
\begin{proof}[Proof of Theorem $\ref{BbmPowers}$]
It suffices to show that $I_{\Delta}$ is a complete intersection ideal 
whenever 
$S/I_{\Delta}^{\ell}$ is Buchsbaum for infinitely many $\ell \ge 1$. 
\par  
By assumption and the above observation, $e(K[\Delta]) \ge 2^c$. 
On the other hand, $S/I_{\Delta}$ is Buchsbaum and thus pure 
by \cite[Theorem 2.6]{HTT}.  
We also have that $\Delta$ is a locally complete intersection complex
by the Goto--Takayama Theorem.  
\par 
Suppose that $\Delta$ is not a complete intersection complex. 
Then by Proposition \ref{LCI-mult}, we have  that $e(K[\Delta]) < 2^c$.
This is a contradiction. 
Hence $\Delta$ must be a complete intersection complex. 
\end{proof}

\begin{exam}  \label{n-gon-nonBbm}
Let $\Delta=\Delta_n$ be the $n$-gon for $n \ge 5$
(or the $n$-pointed path for $n \ge 4$). 
Then $S/I_{\Delta}^{\ell}$ is \textit{not} Buchsbaum for $\ell \ge 6$. 
\end{exam}

\begin{proof}
We consider the case of $n$-gons only. 
Set $I=I_{\Delta} =(X_1X_3,X_1X_4,\ldots,X_{n-2}X_n)$. 
Then $e = e(S/I)=n$, $c= \codim I=n-2$ and $q=\indeg I =2$. 
\par
Suppose that $S/I_{\Delta}^{\ell}$ is Buchsbaum. 
By Corollary \ref{Mult-Bbm}, 
\[
 n=e(S/I) \ge 
\frac{(2\ell+n-4)\cdots (2\ell+1)2\ell(2\ell-1)}{(\ell+n-3)\cdots (\ell+1)\ell}.  
\]
Fix $n \ge 5$ and put $f(\ell)$ to be the right-hand side of the above inequality. 
Then one can easily see that $f(\ell)$ is an increasing function of $\ell$.  
Thus if $\ell \ge 6$, then 
\[
 1 \ge \frac{(n+8)\cdots 12 \cdot 11}{(n+3)\cdots 7\cdot 6} \times \frac{1}{n}
= \frac{(n+8)(n+7)(n+6)(n+5)(n+4)}{10 \cdot 9 \cdot 8  \cdot 7  \cdot 6 \cdot n}. 
\]
Put $g(n)$ to be the right-hand side of the above inequality. 
Then since 
\[
 g(n+1)/g(n) = \frac{n^2+9n}{n^2+5n+4} \ge 1 \quad \text{and} \quad g(5) = 1.02 \cdots > 1   
\]
we get a contradiction. 
\end{proof}

\par 
It is difficult to determine the Buchsbaumness for $S/I^{\ell}$. 

\begin{exam} \label{5-gon2nd}
Let $S=K[X_1,X_2,X_3,X_4,X_5]$ be a polynomial ring. 
Let $I=(X_1X_3,X_1X_4,X_2X_4,X_2X_5,X_3X_5)$ be 
the Stanley--Reisner ideal (of height $3$) of the $5$-gon. 
Then $S/I^2$ is Cohen--Macaulay with $\dim S/I^2 = 2$. 
Indeed, \textbf{Macaulay 2} yields the following minimal free 
resolution of $S/I^2$:
\[
 0 \to S^{10}(-6) \to S^{24}(-5) \to S^{15}(-4) \to S \to S/I^2 \to 0. 
\] 
On the other hand, $\depth S/I^3 = 0$ since 
$X_1X_2X_3X_4X_5 \in I^3 \colon \frm \setminus I^3$. 
We do not know whether $S/I^3$ is Buchsbaum or not. 
\end{exam}
 
\par 
In the following, we give an example of the simplicial complex $\Delta$ for which 
$S/I_{\Delta}^2$ is Buchsbaum 
but \textit{not} Cohen--Macaulay (and this implies that 
$\Delta$ is not a complete intersection complex).   
In order to do that, we use an extension of Hochster's formula 
describing the local cohomology of a monomial ideal; see \cite{Ta}.  
Fix $\ell \ge 1$ and set $G(I_{\Delta}^{\ell}) = \{m_1,\ldots,m_{\mu} \}$. 
Write $m=X_1^{\nu_1(m)} \cdots X_n^{\nu_{n}(m)}$ for any monomial $m$  
in $S=K[X_1,\ldots,X_n]$.  
For a vector ${\bf a} = (a_1,\ldots,a_n) \in \bbZ^n$, we put 
\[
 G_a =\{i \in V \;:\; a_i < 0 \}.
\] 
Then we define the simplicial complex
 $\Delta_{\bf a}(I_{\Delta}^{\ell}) \subseteq \Delta$ by 
\[
 \Delta_{\bf a}(I_{\Delta}^{\ell}) = 
\{L \setminus G_a \;:\: 
G_a \subseteq L \in \Delta, \; \text{$L$ satisfies the condition $(*)$}\},  
\]
where 
\[
(*) \qquad \text{for all $m \in G(I_{\Delta}^{\ell})$, 
there exists an $i \in V \setminus L$  such that $\nu_i(m) > a_i  (\ge 0)$}.
\]
For a graded $S$-module $M$, $F(A,{\bf t}) = \sum_{{\bf a} \in \bbZ^n} \dim_K A_{\bf a} {\bf t}^{\bf a}$ is
called the Hilbert--Poincar\'e series of $M$.  
Then Hochster--Takayama formula (see \cite{Ta}) says that  
\[
 F(H_{\frm}^i(S/I_{\Delta}^{\ell}), {\bf t}) = 
\sum_{F \in \Delta} 
\sum_{\genfrac{}{}{0pt}{}{{\bf a} \in \bbZ^n}{G_a=F, a_i \le \ell-1}} 
\dim_K \widetilde{H}_{i-\sharp(F)-1}(\Delta_{{\bf a}} (I_{\Delta}^{\ell}); K)\; {\bf t}^{\bf a},  
\] 
where $\widetilde{H}_i(\Delta;K)$ denotes the $i$th simplicial reduced homology  of $\Delta$ with 
values in $K$. 
In particular, we have 
\[
 F(H_{\frm}^1(S/I_{\Delta}^{\ell}), {\bf t})
= \sum_{{\bf a} \in \calA} 
\dim_K \widetilde{H}_0 (\Delta_{\bf a} (I_{\Delta}^{\ell});K) {\bf t}^{\bf a} + 
\sum_{i=1}^n \sum_{{\bf a} \in \calA_i} {\bf t}^{\bf a}, 
\]
where 
\begin{eqnarray*}
 \calA &=& \{{\bf a} \in \bbZ^n\;:\; 0 \le a_1,\ldots,a_n \le \ell-1,\; 
\text{$\Delta_{{\bf a}}(I_{\Delta}^{\ell})$ is disconnected }\}; \\
 \calA_i & = & \{{\bf a} \in \bbZ^n\;:\; 0 \le a_1,\ldots, \widehat{a_i} \ldots, a_n \le \ell-1,\; 
\Delta_{{\bf a}}(I_{\Delta}^{\ell}) =\{\emptyset\} \} 
\end{eqnarray*}
for each $i=1,\ldots,n$. 

\par \vspace{3mm}
\begin{exam} \label{BbmnonCM} 
Let $S=K[X_1,X_2,X_3,X_4]$ be a polynomial ring over a field $K$. 
Let $I=(X_1X_3,X_1X_4,X_2X_4)$ be the Stanley--Reisner ideal 
of the $4$-pointed path $\Delta$.  
\par 
Then $S/I^2$ is Buchsbaum but not Cohen--Macaulay. 
In fact, $\dim S/I^2 = 2$, $\depth S/I^2=1$ and $\dim_K H_{\frm}^1(S/I^2)=1$.  
\end{exam}

\begin{proof}
The ideal $I$ can be considered as the edge ideal of some bipartite graph $G$.  
Thus we have $I^2=I^{(2)}$, the second symbolic power of $I$, by 
\cite[Section 5]{SVV}, and so $H_{\frm}^0(S/I^2) =0$. 

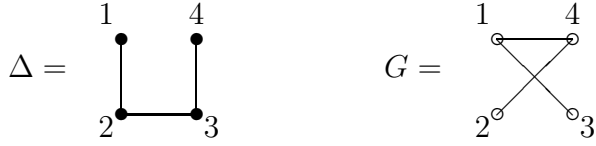
\setlength{\unitlength}{1mm}
\begin{center}
\begin{picture}(100,26)
\put(0,10){$\Delta=$}
\put(15,5){\line(1,0){10}} 
\put(15,5){\line(0,1){10}} 
\put(25,5){\line(0,1){10}} 
\put(14,4){$\bullet$}
\put(24,4){$\bullet$}
\put(14,14){$\bullet$}
\put(24,14){$\bullet$}
\put(12,17){$1$}
\put(12,2){$2$}
\put(26,2){$3$}
\put(24,17){$4$}
\put(50,10){$G=$}
\put(65,15){\line(1,-1){10}} 
\put(65,15){\line(1,0){10}} 
\put(65,5){\line(1,1){10}} 
\put(64,4){$\circ$}
\put(74,4){$\circ$}
\put(64,14){$\circ$}
\put(74,14){$\circ$}
\put(62,17){$1$}
\put(62,2){$2$}
\put(76,2){$3$}
\put(74,17){$4$}
\end{picture}
\end{center}
\par \vspace{2mm}
Hence it suffices to show that  
$\frm H_{\frm}^1(S/I^2)=0$ and $H_{\frm}^1(S/I^2) \ne 0$. 
We first show the following claim. 
Put $\Delta_{{\bf a}} = \Delta_{{\bf a}}(I^2)$ for simplicity. 
 
\begin{description}
\item[Claim 1] $\calA = \{(1,0,0,1)\}$ and $\Delta_{(1,0,0,1)} $ is 
spanned by \{$\{(1,2)\}, \{3,4\}\}$. 
(This implies that $Kt_1t_4 \subseteq H_{\frm}^1(S/I^2)$.)
\end{description}
First of all, we define monomials $m_1,\ldots,m_6$ as follows$:$
\[
\begin{array}{|c|c|c|c|c|c|c|} \hline 
& m_1 & m_2   & m_3 & m_4   & m_5 & m_6 \\ \hline  
\nu_1(m) & 2 & 2 & 2 & 1 & 1 & 0 \\ \hline 
\nu_2(m) & 0 & 0 & 0 & 1 & 1 & 2 \\ \hline 
\nu_3(m) & 2 & 1 & 0 & 1 & 0 & 0 \\ \hline 
\nu_4(m) & 0 & 1 & 2 & 1 & 2 & 2  \\ \hline 
\end{array}
\]
Namely, 
\[
 G(I^2)=\{X_1^2X_3^2,X_1^2X_3X_4,X_1^2X_4^2,X_1X_2X_3X_4,X_1X_2X_4^2,X_2^2X_4^2\}.
\]
Fix ${\bf a}=(a_1,a_2,a_3,a_4) \in (\bbZ \cap \{0,1\})^4$. 
As $\nu_3(m_4) = \nu_4(m_4) =1$, 
it follows that $\{1,2\} \in \Delta_{{\bf a}}$ if and only if 
$a_3 =0$ or $a_4=0$. 
Similarly, $\{3,4\} \in \Delta_{{\bf a}}$ if and only if $a_1=0$ or $a_2=0$. 
If $\sharp\{i \;:\; 1 \le i \le 4,\; a_i =1\} \ge 3$, 
then $\Delta_{{\bf a}}=\emptyset$. 
So, we may assume that $\sharp\{i \;:\; 1 \le i \le 4,\; a_i=1\} \le 2$ 
and $a_1 \ge a_4$.  
\par 
If $\{2,3\} \notin \Delta_{{\bf a}}$, then $a_1=a_4=1$. That is, ${\bf a}=(1,0,0,1)$.  
Indeed, $\Delta_{(1,0,0,1)} = \langle \{1,2\},\,\{3,4\} \rangle$ 
is disconnected.  
Otherwise, $\{2,3\} \in \Delta_{{\bf a}}(I^2)$. Then $(a_1,a_4)=(0,0)$ or $(1,0)$. 
In these cases, we have 
\[
\Delta_{(0,*,*,0)}=\Delta_{(1,0,0,0)}=\Delta_{(1,0,1,0)}=\Delta, \quad
\Delta_{(1,1,0,0)} = \langle \{1,2\},\,\{2,3\} \rangle.  
\]
In particular, $\Delta_{{\bf a}}$ is connected in any case. 
Therefore we proved Claim 1. 
\par \vspace{2mm}
Next, we show the following claim. 
\begin{description}
 \item[Claim 2] $\calA_1 = \calA_2 = \calA_3=\calA_4 = \emptyset$. 
\end{description}
To see $\calA_1 =\emptyset$, let ${{\bf a}} = (a_1,a_2,a_3,a_4) \in \bbZ^4$ 
such that $a_1 < 0$, $0 \le a_2,a_3,a_4 \le 1$. 
Note that 
\[
 \Delta_{{\bf a}}(I^2) =\{L \setminus \{1\} \;:\; \{1\} \subseteq L \in \Delta,\, 
\text{$L$ satisfies $(*)$}\}
\]
and that $\{1\} \subseteq L \in \Delta$ if and only if $L=\{1\}$ or $\{1,2\}$. 
By a similar argument as in the proof of the claim 1, we obtain that 
\[
 \{2\} = \{1,2\} \setminus \{1\} 
\in \Delta_{{\bf a}}(I^2) 
\Longleftrightarrow a_3=0 \;\text{or} \; a_4=0.
\] 
Then $\Delta_{{\bf a}}(I^2)= \{\emptyset, \, \{2\}\} \ne \{\emptyset \}$. 
\par 
Now suppose that $a_3=a_4=1$. 
Then $\emptyset \notin \Delta_{{\bf a}}(I^2)$ because $m_2 =X_1^2X_3X_4 \in G(I^2)$. 
This yields that $\Delta_{{\bf a}}(I^2) \ne \{\emptyset\}$. 
Therefore $\calA_1= \emptyset$. 
Similarly, one has $\calA_2=\calA_3=\calA_4=\emptyset$.   

\par \vspace{2mm}
The above two claims imply that $H_{\frm}^1(S/I^2) \cong Kt_1t_4$, 
as required. 
\end{proof}

\vspace{3mm}
\begin{quest} \label{Indeg-1}
Can you replace Buchsbaumness with 
quasi-Buchsbaumness in Theorem \ref{BbmPowers}?
\end{quest}

\begin{quest} \label{General}
Let $I$ be a generically complete intersection 
homogeneous ideal of a polynomial ring $S$. 
If $S/I^{\ell}$ is Buchsbaum for all $\ell \ge 1$, 
then is $I$ a complete intersection ideal?   
\end{quest}



\end{document}